\documentclass[10pt,a4paper]{amsart}

\usepackage{amsmath}
\usepackage{amssymb}
\usepackage{amsthm}
\usepackage{amsfonts}
\usepackage{microtype}
\usepackage{mathrsfs}
\usepackage{graphicx}
\usepackage{color}
\usepackage{latexsym}
\usepackage{rotating}
\usepackage{xspace}
\usepackage[all]{xy}
\usepackage{longtable}
\usepackage{leftidx}
\usepackage{mathtools}
\usepackage{titlesec}
\usepackage{tikz}
\usetikzlibrary{calc}
\usetikzlibrary{arrows}
\usetikzlibrary{decorations.markings}
\usepackage{geometry}

\newtheorem{theorem}{Theorem}[section]
\numberwithin{equation}{theorem}
\newtheorem{lemma}[theorem]{Lemma}

\newtheorem{corollary}[theorem]{Corollary}

\theoremstyle{definition}
\newtheorem{definition}[theorem]{Definition}

\theoremstyle{conjecture}

\newcommand{\im}{\operatorname{im}}

\newcommand{\Spec}{\operatorname{Spec}}

\newcommand{\rank}{\operatorname{rank}}

\newcommand{\Supp}{\operatorname{Supp}}

\newcommand{\Hom}{\operatorname{Hom}}

\newcommand{\Max}{\operatorname{Max}}

\newcommand{\suchthat}{\;\ifnum\currentgrouptype=16 \middle\fi|\;}

\newenvironment{prf}[1][Proof]{\begin{proof}[\bf #1]}{\end{proof}}

\makeatletter
\newcommand{\holim@}[2]{%
  \vtop{\m@th\ialign{##\cr
    \hfil$#1\operator@font holim$\hfil\cr
    \noalign{\nointerlineskip\kern1.5\ex@}#2\cr
    \noalign{\nointerlineskip\kern-\ex@}\cr}}%
}
\newcommand{\holim}{%
  \mathop{\mathpalette\holim@{\rightarrowfill@\textstyle}}\nmlimits@
}
\makeatother

\makeatletter
\def\@secnumfont{\bfseries}
\def\section{\@startsection{section}{1}%
  \z@{.7\linespacing\@plus\linespacing}{.5\linespacing}%
  {\normalfont\Large\bfseries\filcenter}}
\def\subsection{\@startsection{subsection}{2}%
  \z@{.5\linespacing\@plus.7\linespacing}{-.5em}%
  {\normalfont\large\bfseries}}
\makeatother

\begin{document}

\author[H. Faridian]{Hossein Faridian}

\title[A Note on Projective Modules]
{A Note on Projective Modules}

\address{H. Faridian, School of Mathematical and Statistical Sciences, Clemson University, SC 29634, USA.}
\email{hfaridi@g.clemson.edu, h.faridian@yahoo.com}



\begin{abstract}
This expository note delves into the theory of projective modules parallel to the one developed for injective modules by Matlis. Given a perfect ring $R$, we present a characterization of indecomposable projective $R$-modules and describe a one-to-one correspondence between the projective indecomposable $R$-modules and the simple $R$-modules.
\end{abstract}

\maketitle

\sloppy

\section{Introduction}

Throughout this note, $R$ denotes an associative ring with identity. All modules are assumed to be left and unitary.

There is a well-known theorem in the representation theory of finite-dimensional algebras over a field which provides a 1-1 correspondence between the isomorphism classes of indecomposable projective modules over the algebra on the one hand, and the isomorphism classes of simple modules over the algebra on the other hand.
More specifically, let $k$ be a field and $A$ a finite-dimensional $k$-algebra. Then every simple $A$-module is a quotient of some indecomposable projective
$A$-module which is unique up to isomorphism. Conversely, for every indecomposable projective $A$-module $P$, there is a simple $A$-module, unique up to isomorphism, that is a quotient of $P$ by some maximal submodule; see \cite{Le}. For a generalization of this result to Artin algebras, see \cite[Corollary 4.5]{ARS}, and to perfect rings, see \cite[Proposition 5]{Sh}. The purpose of this note is to present a novel descriptive proof of this theorem in the following more comprehensive form; see Theorems \ref{3.5} and \ref{3.8}.

\begin{theorem} \label{1.1}
Let $R$ be a perfect ring, and $P$ a nonzero $R$-module. Then the following assertions are equivalent:
\begin{enumerate}
\item[(i)] $P$ is an indecomposable projective $R$-module.
\item[(ii)] $P$ is a sum-irreducible projective $R$-module.
\item[(iii)] $P$ is the projective cover of its every nonzero quotient module.
\item[(iv)] $P$ is the projective cover of $R/\mathfrak{m}$ for some maximal left ideal $\mathfrak{m}$ of $R$.
\end{enumerate}
Further if $R$ is commutative, then the above assertions are equivalent to the following one:
\begin{enumerate}
\item[(v)] $P \cong R_{\mathfrak{m}}$ for some maximal ideal $\mathfrak{m}$ of $R$.
\end{enumerate}
\end{theorem}

\begin{theorem} \label{1.2}
Let $R$ be a perfect ring. There is a one-to-one correspondence between the indecomposable projective $R$-modules and the simple $R$-modules.
\end{theorem}

\section{Preliminaries}

We begin with reminding the classical notion of a projective cover.

\begin{definition} \label{2.1}
Let $R$ be a ring, and $M$ an $R$-module. Then:
\begin{enumerate}
\item[(i)] A submodule $N$ of $M$ is said to be \textit{superfluous}, written as $N \subseteq_{sup} M$, if whenever $N + L = M$ for some submodule $L$ of $M$,
then we have $L=M$.
\item[(ii)] By a \textit{projective cover} of $M$, we mean a projective $R$-module $P$ together with an epimorphism $\pi: P \rightarrow M$ such that $\ker \pi
\subseteq_{sup} P$.
\item[(ii)] $R$ is said to be \textit{perfect} if every $R$-module has a projective cover.
\end{enumerate}
\end{definition}

Enochs and Jenda have defined the general notion of a cover as follows; see \cite[Definition 5.1.1]{EJ}.

\begin{definition} \label{2.2}
Let $R$ be a ring, $M$ an $R$-module, and $\mathcal{A}$ a class of $R$-modules. By an \textit{$\mathcal{A}$-cover} of $M$, we mean an $R$-module $A\in \mathcal{A}$
together with an $R$-homomorphism $\varphi: A \rightarrow M$ that satisfy the following conditions:
\begin{enumerate}
\item[(i)] Any $R$-homomorphism $\psi:B \rightarrow M$ with $B \in \mathcal{A}$ factors through $\varphi$, i.e. there is an $R$-homomorphism $f:B\rightarrow A$ that
makes the following diagram commutative:
\[\begin{tikzpicture}[every node/.style={midway}]
  \matrix[column sep={3em}, row sep={3em}]
  {\node(1) {$$}; & \node(2) {$B$}; \\
  \node(3) {$A$}; & \node(4) {$M$};\\};
  \draw[decoration={markings,mark=at position 1 with {\arrow[scale=1.5]{>}}},postaction={decorate},shorten >=0.5pt] (2) -- (4) node[anchor=west] {$\psi$};
  \draw[decoration={markings,mark=at position 1 with {\arrow[scale=1.5]{>}}},postaction={decorate},shorten >=0.5pt] (3) -- (4) node[anchor=south] {$\varphi$};
  \draw[dashed,decoration={markings,mark=at position 1 with {\arrow[scale=1.5]{>}}},postaction={decorate},shorten >=0.5pt] (2) -- (3) node[anchor=south] {$f$};
\end{tikzpicture}\]
\item[(ii)] Any $R$-homomorphism $f:A\rightarrow A$ that makes the diagram
\[\begin{tikzpicture}[every node/.style={midway}]
  \matrix[column sep={3em}, row sep={3em}]
  {\node(1) {$$}; & \node(2) {$A$}; \\
  \node(3) {$A$}; & \node(4) {$M$};\\};
  \draw[decoration={markings,mark=at position 1 with {\arrow[scale=1.5]{>}}},postaction={decorate},shorten >=0.5pt] (2) -- (4) node[anchor=west] {$\varphi$};
  \draw[decoration={markings,mark=at position 1 with {\arrow[scale=1.5]{>}}},postaction={decorate},shorten >=0.5pt] (3) -- (4) node[anchor=south] {$\varphi$};
  \draw[dashed,decoration={markings,mark=at position 1 with {\arrow[scale=1.5]{>}}},postaction={decorate},shorten >=0.5pt] (2) -- (3) node[anchor=south] {$f$};
\end{tikzpicture}\]
commutative is an automorphism.
\end{enumerate}
\end{definition}

The classical definition of a projective cover amounts to the modern definition as recorded in the following result for the sake of bookkeeping.

\begin{lemma} \label{2.3}
Let $R$ be a ring, $M$ an $R$-module, $P$ a projective $R$-module, and $\pi:P \rightarrow M$ an $R$-homomorphism. Let $\mathcal{P}(R)$ denote the class of projective
$R$-modules. Then the following assertions are equivalent:
\begin{enumerate}
\item[(i)] $\pi:P \rightarrow M$ is a projective cover of $M$.
\item[(ii)] $\pi:P \rightarrow M$ is a $\mathcal{P}(R)$-cover of $M$.
\end{enumerate}
\end{lemma}

\begin{prf}
See \cite[Theorem 1.2.12]{Xu}.
\end{prf}

In light of Lemma \ref{2.3}, the projective cover of a given $R$-module $M$ is unique up to isomorphism if it exists, in the sense that if $\pi:P \rightarrow M$
and $\pi^{\prime}: P^{\prime} \rightarrow M$ are two projective covers of $M$, then there is an isomorphism $P^{\prime} \xrightarrow{\cong} P$ that makes the
following diagram commutative:
\[\begin{tikzpicture}[every node/.style={midway}]
  \matrix[column sep={3em}, row sep={3em}]
  {\node(1) {$$}; & \node(2) {$P^{\prime}$}; \\
  \node(3) {$P$}; & \node(4) {$M$};\\};
  \draw[decoration={markings,mark=at position 1 with {\arrow[scale=1.5]{>}}},postaction={decorate},shorten >=0.5pt] (2) -- (4) node[anchor=west] {$\pi^{\prime}$};
  \draw[decoration={markings,mark=at position 1 with {\arrow[scale=1.5]{>}}},postaction={decorate},shorten >=0.5pt] (3) -- (4) node[anchor=south] {$\pi$};
  \draw[dashed,decoration={markings,mark=at position 1 with {\arrow[scale=1.5]{>}}},postaction={decorate},shorten >=0.5pt] (2) -- (3) node[anchor=south] {$\cong$};
\end{tikzpicture}\]
Accordingly, we denote a choice of a projective cover of $M$ by $\pi_{M}: P_{R}(M) \rightarrow M$ whenever one exists.

\begin{definition} \label{2.4}
Let $R$ be a ring. Then:
\begin{enumerate}
\item[(i)]  A nonzero $R$-module $M$ is said to be \textit{indecomposable} if whenever $M=N_{1}\oplus N_{2}$ for some submodules $N_{1}$ and $N_{2}$ of $M$,
then we have either $N_{1}=0$ or $N_{2}=0$.
\item[(ii)] An $R$-module $M$ is said to be \textit{sum-irreducible} if whenever $M=N_{1} + N_{2}$ for some submodules $N_{1}$ and $N_{2}$ of $M$, then we have
either $N_{1}=M$ or $N_{2}=M$.
\end{enumerate}
\end{definition}

By \cite[Theorem 28.4]{AF}, one has the following lemma:

\begin{lemma} \label{2.5}
Let $R$ be a ring, and $\mathfrak{J}(R)$ denote its Jacobson radical. Then the following assertions are equivalent:
\begin{enumerate}
\item[(i)] $R$ is perfect.
\item[(ii)] Every descending chain of principal right ideals of $R$ stabilizes.
\item[(iii)] Every flat $R$-module is projective.
\item[(iv)] $R/\mathfrak{J}(R)$ is semisimple and $\mathfrak{J}(R)$ is left T-nilpotent.
\end{enumerate}
\end{lemma}

Recall that a subset $A$ of a ring $R$ is said to be left T-nilpotent if for every sequence $a_1,a_2,... \in A$, there is a natural number $n$
such that $a_1a_2\cdots a_n=0$.

\section{Main Results}

In this section, we delve into the structure of indecomposable projective $R$-modules.

\begin{lemma} \label{3.1}
Let $R$ be a perfect ring, and $P$ a nonzero projective $R$-module. Then the following assertions are equivalent:
\begin{enumerate}
\item[(i)] $P$ is indecomposable.
\item[(ii)] Given any proper submodule $N$ of $P$, there is an isomorphism $P \xrightarrow{\cong} P_{R}(P/N)$ that makes the following diagram commutative:
\[\begin{tikzpicture}[every node/.style={midway}]
  \matrix[column sep={3em}, row sep={3em}]
  {\node(1) {$$}; & \node(2) {$P_{R}(P/N)$}; \\
  \node(3) {$P$}; & \node(4) {$P/N$};\\};
  \draw[decoration={markings,mark=at position 1 with {\arrow[scale=1.5]{>}}},postaction={decorate},shorten >=0.5pt] (2) -- (4) node[anchor=west] {$\pi_{P/N}$};
  \draw[decoration={markings,mark=at position 1 with {\arrow[scale=1.5]{>}}},postaction={decorate},shorten >=0.5pt] (3) -- (4) node[anchor=south] {$\pi$};
  \draw[dashed,decoration={markings,mark=at position 1 with {\arrow[scale=1.5]{>}}},postaction={decorate},shorten >=0.5pt] (3) -- (2) node[anchor=east] {$\cong$};
\end{tikzpicture}\] where $\pi$ is the natural epimorphism.
\item[(iii)] $P$ is sum-irreducible.
\end{enumerate}
\end{lemma}

\begin{prf}
(i) $\Rightarrow$ (ii): Let $N$ be a proper submodule of $P$. By Lemma \ref{2.3}, the $R$-homomorphism $\pi_{P/N}: P_{R}(P/N) \rightarrow P/N$ is a
$\mathcal{P}(R)$-cover, so we can find an $R$-homomorphism $f:P \rightarrow P_{R}(P/N)$, and also using the projectivity of $P_{R}(P/N)$, we can find
an $R$-homomorphism $g: P_{R}(P/N) \rightarrow P$ that make the following diagram commutative:
\[\begin{tikzpicture}[every node/.style={midway}]
  \matrix[column sep={3em}, row sep={3em}]
  {\node(1) {$$}; & \node(2) {$P_{R}(P/N)$}; & \node(3) {$$}; \\
  \node(4) {$P$}; & \node(5) {$P/N$}; & \node(6) {$0$};\\};
  \draw[decoration={markings,mark=at position 1 with {\arrow[scale=1.5]{>}}},postaction={decorate},shorten >=0.5pt] (2) -- (5) node[anchor=west] {$\pi_{P/N}$};
  \draw[decoration={markings,mark=at position 1 with {\arrow[scale=1.5]{>}}},postaction={decorate},shorten >=0.5pt] (4) -- (5) node[anchor=south] {$\pi$};
  \draw[decoration={markings,mark=at position 1 with {\arrow[scale=1.5]{>}}},postaction={decorate},shorten >=0.5pt] (5) -- (6) node[anchor=south] {$$};
  \draw[transform canvas={xshift=-0.8ex,yshift=0.3ex},dashed,decoration={markings,mark=at position 1 with {\arrow[scale=1.5]{>}}},postaction={decorate},
  shorten >=0.5pt] (2) -- (4) node[anchor=east] {$g$};
  \draw[transform canvas={xshift=0.4ex,yshift=-0.2ex},dashed,decoration={markings,mark=at position 1 with {\arrow[scale=1.5]{>}}},postaction={decorate},
  shorten >=0.5pt] (4) -- (2) node[anchor=west] {$f$};
\end{tikzpicture}\]
We thus have $\pi_{P/N}fg=\pi g=\pi_{P/N}$. Another use of Lemma \ref{2.3} implies that $fg$ is an automorphism, so $f$ is surjective. The short exact sequence
$$0 \rightarrow \ker f \rightarrow P \xrightarrow{f} P_{R}(P/N) \rightarrow 0$$
splits, so $P \cong \ker f \oplus P_{R}(P/N)$. But $P_{R}(P/N) \neq 0$, so the hypothesis implies that $\ker f=0$, whence $f:P \rightarrow P_{R}(P/N)$ is an
isomorphism.

(ii) $\Rightarrow$ (iii): Suppose that $N_{1}$ and $N_{2}$ are two submodules of $P$ such that $P= N_{1}+N_{2}$. If $N_{1}$ is proper, then the hypothesis
implies that there is an isomorphism $f:P \rightarrow P_{R}(P/N_{1})$ that makes the following diagram commutative:
\[\begin{tikzpicture}[every node/.style={midway}]
  \matrix[column sep={3em}, row sep={3em}]
  {\node(1) {$$}; & \node(2) {$P_{R}(P/N_{1})$}; \\
  \node(3) {$P$}; & \node(4) {$P/N_{1}$};\\};
  \draw[decoration={markings,mark=at position 1 with {\arrow[scale=1.5]{>}}},postaction={decorate},shorten >=0.5pt] (2) -- (4) node[anchor=west] {$\pi_{P/N_{1}}$};
  \draw[decoration={markings,mark=at position 1 with {\arrow[scale=1.5]{>}}},postaction={decorate},shorten >=0.5pt] (3) -- (4) node[anchor=south] {$\pi$};
  \draw[dashed,decoration={markings,mark=at position 1 with {\arrow[scale=1.5]{>}}},postaction={decorate},shorten >=0.5pt] (3) -- (2) node[anchor=east] {$f$};
\end{tikzpicture}\]
We thus obtain the following commutative diagram with exact rows:
\[\begin{tikzpicture}[every node/.style={midway}]
  \matrix[column sep={3em}, row sep={3em}]
  {\node(1) {$0$}; & \node(2) {$N_{1}$}; & \node(3) {$P$}; & \node(4) {$P/N_{1}$}; & \node(5) {$0$};\\
  \node(6) {$0$}; & \node(7) {$\ker \pi_{P/N_{1}}$}; & \node(8) {$P_{R}(P/N_{1})$}; & \node(9) {$P/N_{1}$}; & \node(10) {$0$};\\};
  \draw[decoration={markings,mark=at position 1 with {\arrow[scale=1.5]{>}}},postaction={decorate},shorten >=0.5pt] (1) -- (2) node[anchor=south] {$$};
  \draw[decoration={markings,mark=at position 1 with {\arrow[scale=1.5]{>}}},postaction={decorate},shorten >=0.5pt] (2) -- (3) node[anchor=south] {$$};
  \draw[decoration={markings,mark=at position 1 with {\arrow[scale=1.5]{>}}},postaction={decorate},shorten >=0.5pt] (3) -- (4) node[anchor=south] {$\pi$};
  \draw[decoration={markings,mark=at position 1 with {\arrow[scale=1.5]{>}}},postaction={decorate},shorten >=0.5pt] (4) -- (5) node[anchor=south] {$$};
  \draw[decoration={markings,mark=at position 1 with {\arrow[scale=1.5]{>}}},postaction={decorate},shorten >=0.5pt] (6) -- (7) node[anchor=south] {$$};
  \draw[decoration={markings,mark=at position 1 with {\arrow[scale=1.5]{>}}},postaction={decorate},shorten >=0.5pt] (7) -- (8) node[anchor=south] {$$};
  \draw[decoration={markings,mark=at position 1 with {\arrow[scale=1.5]{>}}},postaction={decorate},shorten >=0.5pt] (8) -- (9) node[anchor=south] {$\pi_{P/N_{1}}$};
  \draw[decoration={markings,mark=at position 1 with {\arrow[scale=1.5]{>}}},postaction={decorate},shorten >=0.5pt] (9) -- (10) node[anchor=south] {$$};
  \draw[dashed,decoration={markings,mark=at position 1 with {\arrow[scale=1.5]{>}}},postaction={decorate},shorten >=0.5pt] (2) -- (7) node[anchor=west] {$\cong$};
  \draw[decoration={markings,mark=at position 1 with {\arrow[scale=1.5]{>}}},postaction={decorate},shorten >=0.5pt] (3) -- (8) node[anchor=west] {$f$};
  \draw[double distance=1.5pt] (4) -- (9) node[anchor=west] {};
\end{tikzpicture}\]
Since $\ker \pi_{P/N_{1}} \subseteq_{sup} P_{R}(P/N_{1})$, we deduces that $N_{1} \subseteq _{sup} P$. It follows that $N_{2}=P$, so $P$ is sum-irreducible.

(iii) $\Rightarrow$ (i): Suppose that $P=N_{1} \oplus N_{2}$ for some submodules $N_{1}$ and $N_{2}$ of $P$. The hypothesis implies that either $N_{1}=P$ or
$N_{2}=P$. It follows that either $N_{2}=0$ or $N_{1}=0$, so $P$ is indecomposable.
\end{prf}

\begin{corollary} \label{3.2}
Let $R$ be a perfect ring, and $M$ an $R$-module. Then $M$ is sum-irreducible if and only if $P_{R}(M)$ is indecomposable.
\end{corollary}

\begin{prf}
Let $\pi_{M}: P_{R}(M) \rightarrow M$ be a projective cover of $M$. Suppose that $M$ is sum-irreducible, and let $P_{R}(M)= P_{1} + P_{2}$ for some submodules
$P_{1}$ and $P_{2}$ of $P_{R}(M)$. Then
$$M= \pi_{M}\left(P_{R}(M)\right)= \pi_{M}(P_{1}) + \pi_{M}(P_{2}).$$
Without loss of generality, we can conclude that $M=\pi_{M}(P_{1})$. It follows that $P_{R}(M)= P_{1} + \ker \pi_{M}$. But $\ker \pi_{M} \subseteq_{sup} P_{R}(M)$,
so we infer that $P_{R}(M)=P_{1}$. It follows that $P_{R}(M)$ is sum-irreducible. Therefore, by Lemma \ref{3.1}, $P_{R}(M)$ is indecomposable.

Conversely, suppose that $P_{R}(M)$ is indecomposable. Hence by Lemma \ref{3.1}, $P_{R}(M)$ is sum-irreducible. Assume that $M= N_{1} + N_{2}$. It is easy to
see that $$P_{R}(M)=\pi^{-1}_{M}(N_{1})+\pi^{-1}_{M}(N_{2})+\ker \pi_{M}.$$
But $\ker \pi_{M} \subseteq_{sup} P_{R}(M)$, so we infer that $P_{R}(M)= \pi^{-1}_{M}(N_{1}) + \pi^{-1}_{M}(N_{2})$. Without loss of generality, we can conclude that
$P_{R}(M)= \pi^{-1}_{M}(N_{1})$. This shows that $M=N_{1}$, so $M$ is sum-irreducible.
\end{prf}

\begin{corollary} \label{3.3}
Let $R$ be a perfect ring, and $P$ a nonzero projective $R$-module. Then $P$ is indecomposable if and only if $P \cong P_{R}(S)$ for some simple $R$-module $S$.
\end{corollary}

\begin{prf}
Suppose that $P$ is indecomposable. By \cite[Proposition 17.14]{AF}, $P$ has a maximal submodule $N$. Then by Lemma \ref{3.1}, $P \cong P_{R}(P/N)$. It is
clear that $S:= P/N$ is a simple $R$-module.

Conversely, suppose that $P \cong P_{R}(S)$ for some simple $R$-module $S$. Obviously, $S$ is sum-irreducible. Now, Corollary \ref{3.2} implies that $P_{R}(S)$ is
indecomposable.
\end{prf}

\begin{lemma} \label{3.4}
Let $R$ be a commutative perfect ring. Then the following assertions hold:
\begin{enumerate}
\item[(i)] $\Spec(R)=\Max(R)$.
\item[(ii)] $\Supp_{R}(R_{\mathfrak{m}})= \{\mathfrak{m}\}$ for every $\mathfrak{m} \in \Max(R)$.
\item[(iii)] $R_{\mathfrak{m}}$ is an indecomposable $R$-module for every $\mathfrak{m} \in \Max(R)$.
\end{enumerate}
\end{lemma}

\begin{prf}
(i): Let $\mathfrak{p} \in \Spec(R)$, and $0 \neq a+ \mathfrak{p} \in R/\mathfrak{p}$. By Lemma \ref{2.5}, the descending chain $(a) \supseteq (a^{2}) \supseteq
\cdots$ of principal ideals of $R$ stabilizes, i.e. there is an integer $n \geq 1$ such that $(a^{n})=(a^{n+1})=\cdots$. In particular, $a^{n} = ra^{n+1}$ for
some $r\in R$. It follows that $(1-ra)a^{n}=0$. As $a \notin \mathfrak{p}$, we get $1-ra \in \mathfrak{p}$, which yields that $(r+\mathfrak{p})(a+\mathfrak{p})
=1+\mathfrak{p}$, i.e. $a+\mathfrak{p} \in (R/\mathfrak{p})^{\times}$. Hence $R/\mathfrak{p}$ is a field, so $\mathfrak{p} \in \Max(R)$.

(ii): Let $\mathfrak{m}, \mathfrak{n} \in \Max(R)$ be such that $\mathfrak{m} \neq \mathfrak{n}$. Set $$\mathfrak{J}_{\mathfrak{u}}(R):= \bigcap _{\mathfrak{v}
\in \Max(R) \backslash \{\mathfrak{u}\}} \mathfrak{v},$$ for any $\mathfrak{u}\in \Max(R)$. By Lemma \ref{2.5}, the ring $R/\mathfrak{J}(R)$ is semisimple. This
implies that $R$ is semilocal, and so for any $i \geq 1$, there are elements $a_{i}\in \mathfrak{J}_{\mathfrak{m}}(R)^{i} \backslash \mathfrak{m}$
and $b_{i}\in \mathfrak{J}_{\mathfrak{n}}(R)^{i} \backslash \mathfrak{n}$. Set $c_{i}:=a_{i}b_{i}$ for every $i \geq 1$. It is clear that $c_{i}\in \mathfrak{J}(R)$
for every $i \geq 1$. By Lemma \ref{2.5}, $\mathfrak{J}(R)$ is T-nilpotent. It follows that there is an integer $n \geq 1$ such that $c_{1}c_{2}\cdots c_{n}=0$.
Set $a:=a_{1}a_{2}\cdots a_{n}$ and $b:=b_{1}b_{2}\cdots b_{n}$. Then it is obvious that $a \notin \mathfrak{m}$, $b \notin \mathfrak{n}$, and $ab=0$. This shows
that $(R_{\mathfrak{m}})_{\mathfrak{n}}=0$, so $\Supp_{R}(R_{\mathfrak{m}})= \{\mathfrak{m}\}$.

(iii): Suppose that $R_{\mathfrak{m}}=N \oplus N^{\prime}$ for some $R$-submodules $N$ and $N^{\prime}$ of $R_{\mathfrak{m}}$. Localizing at $\mathfrak{m}$, we
get $R_{\mathfrak{m}}=N_{\mathfrak{m}}\oplus N^{\prime}_{\mathfrak{m}}$. But since $R_{\mathfrak{m}}$ is a local ring, it follows that $R_{\mathfrak{m}}$ is
an indecomposable $R_{\mathfrak{m}}$-module. Hence, $N_{\mathfrak{m}}=0$ or $N^{\prime}_{\mathfrak{m}}=0$. Say $N_{\frak m}=0$. But, (ii) implies that
$\Supp_{R}(N)\subseteq \{\mathfrak{m}\}$. Therefore,  $N=0$.
\end{prf}

\begin{theorem} \label{3.5}
Let $R$ be a perfect ring, and $P$ a nonzero $R$-module. Then the following assertions are equivalent:
\begin{enumerate}
\item[(i)] $P$ is an indecomposable projective $R$-module.
\item[(ii)] $P$ is a sum-irreducible projective $R$-module.
\item[(iii)] $P$ is the projective cover of its every nonzero quotient module.
\item[(iv)] $P \cong P_{R}(R/\mathfrak{m})$ for some maximal left ideal $\mathfrak{m}$ of $R$.
\end{enumerate}
Further if $R$ is commutative, then the above assertions are equivalent to the following one:
\begin{enumerate}
\item[(v)] $P \cong R_{\mathfrak{m}}$ for some maximal ideal $\mathfrak{m}$ of $R$.
\end{enumerate}
\end{theorem}

\begin{prf} Lemma \ref{3.1} yields the equivalence of (i), (ii) and (iii).

(i) $\Rightarrow$ (iv): Corollary \ref{3.3} implies that $P \cong P_{R}(S)$ for some simple $R$-module $S$. But then $S \cong R/\mathfrak{m}$ for some maximal left
ideal $\mathfrak{m}$ of $R$.

(iv) $\Rightarrow$ (v): Suppose that $P \cong P_{R}(R/\mathfrak{m})$ for some maximal ideal $\mathfrak{m}$ of $R$. By Lemma \ref{3.4}, $R_{\mathfrak{m}}$ is an
indecomposable $R$-module. Moreover, Lemma \ref{2.5} yields that $R_{\mathfrak{m}}$ is a projective $R$-module. By Corollary \ref{3.3}, we have $R_{\mathfrak{m}}
\cong P_{R}(R/\mathfrak{n})$ for some maximal ideal $\mathfrak{n}$ of $R$. If $\mathfrak{n} \neq \mathfrak{m}$, then take any $a \in \mathfrak{n} \backslash
\mathfrak{m}$. Therefore, the $R$-homomorphism $a1^{R_{m}}: R_{\mathfrak{m}} \rightarrow R_{\mathfrak{m}}$ is an isomorphism. The commutative diagram
\[\begin{tikzpicture}[every node/.style={midway}]
  \matrix[column sep={3em}, row sep={3em}]
  {\node(1) {$R_{\mathfrak{m}}$}; & \node(2) {$R_{\mathfrak{m}}$};\\
  \node(3) {$P_{R}(R/\mathfrak{n})$}; & \node (4) {$P_{R}(R/\mathfrak{n})$};\\
  \node(5) {$R/\mathfrak{n}$}; & \node (6) {$R/\mathfrak{n}$};\\};
  \draw[decoration={markings,mark=at position 1 with {\arrow[scale=1.5]{>}}},postaction={decorate},shorten >=0.5pt] (1) -- (2) node[anchor=south] {$a$};
  \draw[decoration={markings,mark=at position 1 with {\arrow[scale=1.5]{>}}},postaction={decorate},shorten >=0.5pt] (1) -- (3) node[anchor=west] {$\cong$};
  \draw[decoration={markings,mark=at position 1 with {\arrow[scale=1.5]{>}}},postaction={decorate},shorten >=0.5pt] (2) -- (4) node[anchor=west] {$\cong$};
  \draw[decoration={markings,mark=at position 1 with {\arrow[scale=1.5]{>}}},postaction={decorate},shorten >=0.5pt] (3) -- (4) node[anchor=south] {$a$};
  \draw[decoration={markings,mark=at position 1 with {\arrow[scale=1.5]{>}}},postaction={decorate},shorten >=0.5pt] (3) -- (5) node[anchor=west]
  {$\pi_{R/\mathfrak{n}}$};  \draw[decoration={markings,mark=at position 1 with {\arrow[scale=1.5]{>}}},postaction={decorate},shorten >=0.5pt] (4) -- (6)
  node[anchor=west] {$\pi_{R/\mathfrak{n}}$};
  \draw[decoration={markings,mark=at position 1 with {\arrow[scale=1.5]{>}}},postaction={decorate},shorten >=0.5pt] (5) -- (6) node[anchor=south] {$a$};
\end{tikzpicture}\]
shows that the $R$-homomorphism $a1^{R/\mathfrak{n}}:R/\mathfrak{n} \rightarrow R/\mathfrak{n}$ is surjective, which is a contradiction. Hence,
$\mathfrak{n}=\mathfrak{m}$.

(v) $\Rightarrow$ (i): Follows from Lemmas \ref{3.4} and \ref{2.5}.
\end{prf}

The following result may be considered as dual to Matlis Theorem which asserts that over a commutative noetherian ring $R$, every injective $R$-module decomposes
uniquely into a direct sum of indecomposable injective $R$-modules; see \cite{Ma}.

\begin{corollary} \label{3.6}
Let $R$ be a commutative artinian ring, and $P$ a projective $R$-module. Then $$P \cong \prod_{\mathfrak{p}\in \Spec(R)} P_{R}(R/\mathfrak{p})^{\pi_{0}(
\mathfrak{p},P)},$$ where $\pi_{0}(\mathfrak{p},P)$ is the zeroth dual Bass number of $P$ with respect to $\mathfrak{p}$, i.e.
$$\pi_{0}(\mathfrak{p},P)=\rank_{R_{\frak p}/{\frak p}R_{\frak p}}\left(R_{\frak p}/{\frak p}R_{\frak p}\otimes_{R_{\frak p}}\Hom_R(R_{\frak p},P)\right).$$
\end{corollary}

\begin{prf}
By \cite[Theorem 1.2.13]{Xu}, $R$ is perfect, and thus \cite[Proposition 3.3.1]{Xu} implies that $P$ is a flat cotorsion $R$-module. It follows from \cite[Theorem 4.1.15]{Xu}
that
$$P \cong \prod_{\mathfrak{p}\in \Spec(R)} T_{\mathfrak{p}},$$
where $T_{\mathfrak{p}}$ is the $\mathfrak{p}R_{\mathfrak{p}}$-adic completion of a free $R_{\mathfrak{p}}$-module of rank $\pi_{0}(\mathfrak{p},P)$. But $\mathfrak{J}(R_{\mathfrak{p}})=\mathfrak{p}R_{\mathfrak{p}}$ is nilpotent, so every $R_{\mathfrak{p}}$-module is $\mathfrak{p}R_{\mathfrak{p}}$-adically complete. It
follows that $T_{\mathfrak{p}}$ is a free $R_{\mathfrak{p}}$-module of rank $\pi_{0}(\mathfrak{p},P)$. Now, the result follows from Theorem \ref{3.5}.
\end{prf}

\begin{lemma} \label{3.7}
Let $R$ be a perfect ring, and $f:M \rightarrow N$ an $R$-homomorphism. Then the following assertions hold:
\begin{enumerate}
\item[(i)] There is an $R$-homomorphism $\tilde{f}: P_{R}(M) \rightarrow P_{R}(N)$ that makes the following diagram commutative:
\[\begin{tikzpicture}[every node/.style={midway}]
  \matrix[column sep={3em}, row sep={3em}]
  {\node(1) {$P_{R}(M)$}; & \node(2) {$P_{R}(N)$}; \\
  \node(3) {$M$}; & \node(4) {$N$};\\};
  \draw[dashed,decoration={markings,mark=at position 1 with {\arrow[scale=1.5]{>}}},postaction={decorate},shorten >=0.5pt] (1) -- (2) node[anchor=south] {$\tilde{f}$};
  \draw[decoration={markings,mark=at position 1 with {\arrow[scale=1.5]{>}}},postaction={decorate},shorten >=0.5pt] (1) -- (3) node[anchor=west] {$\pi_{M}$};
  \draw[decoration={markings,mark=at position 1 with {\arrow[scale=1.5]{>}}},postaction={decorate},shorten >=0.5pt] (2) -- (4) node[anchor=west] {$\pi_{N}$};
  \draw[decoration={markings,mark=at position 1 with {\arrow[scale=1.5]{>}}},postaction={decorate},shorten >=0.5pt] (3) -- (4) node[anchor=south] {$f$};
\end{tikzpicture}\]
\item[(ii)] If $f$ is an epimorphism, then any $R$-homomorphism $g: P_{R}(M) \rightarrow P_{R}(N)$ that makes the diagram in (i) commutative is an epimorphism.\\
\item[(iii)] If $f$ is an isomorphism, then any $R$-homomorphism $g: P_{R}(M) \rightarrow P_{R}(N)$ that makes the diagram in (i) commutative is an isomorphism.
\end{enumerate}
\end{lemma}

\begin{prf}
(i): The existence of $\tilde{f}$ follows readily from the projectivity of $P_{R}(M)$ in light of the following diagram:
\[\begin{tikzpicture}[every node/.style={midway}]
  \matrix[column sep={2em}, row sep={2em}]
  {\node(1) {$$}; & \node(2) {$P_{R}(M)$}; & \node(3) {$$};\\
  \node(4) {$$}; & \node(5) {$M$}; & \node(6) {$$};\\
  \node(7) {$P_{R}(N)$}; & \node(8) {$N$}; & \node(9) {$0$};\\};
  \draw[decoration={markings,mark=at position 1 with {\arrow[scale=1.5]{>}}},postaction={decorate},shorten >=0.5pt] (2) -- (5) node[anchor=west] {$\pi_{M}$};
  \draw[decoration={markings,mark=at position 1 with {\arrow[scale=1.5]{>}}},postaction={decorate},shorten >=0.5pt] (5) -- (8) node[anchor=west] {$f$};
  \draw[decoration={markings,mark=at position 1 with {\arrow[scale=1.5]{>}}},postaction={decorate},shorten >=0.5pt] (7) -- (8) node[anchor=south] {$\pi_{N}$};
  \draw[decoration={markings,mark=at position 1 with {\arrow[scale=1.5]{>}}},postaction={decorate},shorten >=0.5pt] (8) -- (9) node[anchor=south] {$$};
\draw[dashed,decoration={markings,mark=at position 1 with {\arrow[scale=1.5]{>}}},postaction={decorate},shorten >=0.5pt] (2) -- (7) node[anchor=east] {$\tilde{f}$};
\end{tikzpicture}\]

(ii): Suppose that $f$ is an epimorphism and an $R$-homomorphism $g: P_{R}(M) \rightarrow P_{R}(N)$ makes the following diagram commutative:
\[\begin{tikzpicture}[every node/.style={midway}]
  \matrix[column sep={3em}, row sep={3em}]
  {\node(1) {$P_{R}(M)$}; & \node(2) {$P_{R}(N)$}; \\
  \node(3) {$M$}; & \node (4) {$N$};\\};
  \draw[decoration={markings,mark=at position 1 with {\arrow[scale=1.5]{>}}},postaction={decorate},shorten >=0.5pt] (1) -- (2) node[anchor=south] {$g$};
  \draw[decoration={markings,mark=at position 1 with {\arrow[scale=1.5]{>}}},postaction={decorate},shorten >=0.5pt] (1) -- (3) node[anchor=west] {$\pi_{M}$};
  \draw[decoration={markings,mark=at position 1 with {\arrow[scale=1.5]{>}}},postaction={decorate},shorten >=0.5pt] (2) -- (4) node[anchor=west] {$\pi_{N}$};
  \draw[decoration={markings,mark=at position 1 with {\arrow[scale=1.5]{>}}},postaction={decorate},shorten >=0.5pt] (3) -- (4) node[anchor=south] {$f$};
\end{tikzpicture}\]
The diagram shows that $\pi_{N} g$ is surjective. Therefore, it can be seen by inspection that $$\im g + \ker \pi_{N} = P_{R}(N).$$
But $\ker \pi_{N} \subseteq_{sup} P_{R}(N)$, so $\im g = P_{R}(N)$, i.e. $g$ is surjective.

(iii): Suppose that $f$ is an isomorphism and an $R$-homomorphism $g: P_{R}(M) \rightarrow P_{R}(N)$ makes the following diagram commutative:
\[\begin{tikzpicture}[every node/.style={midway}]
  \matrix[column sep={3em}, row sep={3em}]
  {\node(1) {$P_{R}(M)$}; & \node(2) {$P_{R}(N)$}; \\
  \node(3) {$M$}; & \node (4) {$N$};\\};
  \draw[decoration={markings,mark=at position 1 with {\arrow[scale=1.5]{>}}},postaction={decorate},shorten >=0.5pt] (1) -- (2) node[anchor=south] {$g$};
  \draw[decoration={markings,mark=at position 1 with {\arrow[scale=1.5]{>}}},postaction={decorate},shorten >=0.5pt] (1) -- (3) node[anchor=west] {$\pi_{M}$};
  \draw[decoration={markings,mark=at position 1 with {\arrow[scale=1.5]{>}}},postaction={decorate},shorten >=0.5pt] (2) -- (4) node[anchor=west] {$\pi_{N}$};
  \draw[decoration={markings,mark=at position 1 with {\arrow[scale=1.5]{>}}},postaction={decorate},shorten >=0.5pt] (3) -- (4) node[anchor=south] {$f$};
\end{tikzpicture}\]
By (ii), $g$ is surjective. The short exact sequence
$$0 \rightarrow \ker g \rightarrow P_{R}(M) \xrightarrow{g} P_{R}(N) \rightarrow 0$$ splits, so $P_{R}(M)=\ker g +\im g^{\prime}$,
where $g^{\prime}: P_{R}(N) \rightarrow P_{R}(M)$ is an $R$-homomorphism such that $gg^{\prime}=1^{P_{R}(N)}$.
It is clear that $\ker g \subseteq \ker \pi_{M}$, so $P_{R}(M)= \ker \pi_{M} + \im g^{\prime}$. But $\ker \pi_{M} \subseteq_{sup} P_{R}(M)$,
so $P_{R}(M)=\im g^{\prime}$. It follows that $g^{\prime}$ is an isomorphism, so $g$ is an isomorphism.
\end{prf}

The next result extends the main result of \cite{Le} from finite-dimensional algebras over a field to perfect rings.

\begin{theorem} \label{3.8}
Let $R$ be a perfect ring, and $\mathfrak{m},\mathfrak{n}$ two maximal left ideals of $R$. Then $P_{R}(R/\mathfrak{m}) \cong P_{R}(R/\mathfrak{n})$ if and only
if $\mathfrak{m}=\mathfrak{n}$. Hence, there is a one-to-one correspondence between the indecomposable projective $R$-modules and the simple $R$-modules.
\end{theorem}

\begin{prf}
Let $\varphi: P_{R}(R/\mathfrak{m}) \rightarrow P_{R}(R/\mathfrak{n})$ be any isomorphism. Suppose to the contrary that $\mathfrak{m}\neq \mathfrak{n}$, and take
any element $a \in \mathfrak{m} \backslash \mathfrak{n}$. It follows that the $R$-homomorphism $a1^{R/\mathfrak{n}}:R/\mathfrak{n} \rightarrow R/\mathfrak{n}$ is
an isomorphism. In view of Lemma \ref{3.7}, the commutative diagram
\[\begin{tikzpicture}[every node/.style={midway}]
  \matrix[column sep={3em}, row sep={3em}]
  {\node(1) {$P_{R}(R/\mathfrak{n})$}; & \node(2) {$P_{R}(R/\mathfrak{n})$}; \\
  \node(3) {$R/\mathfrak{n}$}; & \node (4) {$R/\mathfrak{n}$};\\};
  \draw[decoration={markings,mark=at position 1 with {\arrow[scale=1.5]{>}}},postaction={decorate},shorten >=0.5pt] (1) -- (2) node[anchor=south] {$a$};
  \draw[decoration={markings,mark=at position 1 with {\arrow[scale=1.5]{>}}},postaction={decorate},shorten >=0.5pt] (1) -- (3) node[anchor=west] {$\pi_{R/\mathfrak{n}}$};
  \draw[decoration={markings,mark=at position 1 with {\arrow[scale=1.5]{>}}},postaction={decorate},shorten >=0.5pt] (2) -- (4) node[anchor=west] {$\pi_{R/\mathfrak{n}}$};
  \draw[decoration={markings,mark=at position 1 with {\arrow[scale=1.5]{>}}},postaction={decorate},shorten >=0.5pt] (3) -- (4) node[anchor=south] {$a$};
\end{tikzpicture}\]
shows that the $R$-homomorphism $a1^{P_{R}(R/\mathfrak{n})}:P_{R}(R/\mathfrak{n}) \rightarrow P_{R}(R/\mathfrak{n})$ is an isomorphism. Therefore, the commutative diagram
\[\begin{tikzpicture}[every node/.style={midway}]
  \matrix[column sep={3em}, row sep={3em}]
  {\node(1) {$P_{R}(R/\mathfrak{m})$}; & \node(2) {$P_{R}(R/\mathfrak{m})$}; \\
  \node(3) {$P_{R}(R/\mathfrak{n})$}; & \node (4) {$P_{R}(R/\mathfrak{n})$};\\};
  \draw[decoration={markings,mark=at position 1 with {\arrow[scale=1.5]{>}}},postaction={decorate},shorten >=0.5pt] (1) -- (2) node[anchor=south] {$a$};
  \draw[decoration={markings,mark=at position 1 with {\arrow[scale=1.5]{>}}},postaction={decorate},shorten >=0.5pt] (1) -- (3) node[anchor=west] {$\varphi$};
  \draw[decoration={markings,mark=at position 1 with {\arrow[scale=1.5]{>}}},postaction={decorate},shorten >=0.5pt] (2) -- (4) node[anchor=west] {$\varphi$};
  \draw[decoration={markings,mark=at position 1 with {\arrow[scale=1.5]{>}}},postaction={decorate},shorten >=0.5pt] (3) -- (4) node[anchor=south] {$a$};
\end{tikzpicture}\]
yields that the $R$-homomorphism $a1^{P_{R}(R/\mathfrak{m})}:P_{R}(R/\mathfrak{m}) \rightarrow P_{R}(R/\mathfrak{m})$ is an isomorphism. Now, the commutative diagram
\[\begin{tikzpicture}[every node/.style={midway}]
  \matrix[column sep={3em}, row sep={3em}]
  {\node(1) {$P_{R}(R/\mathfrak{m})$}; & \node(2) {$P_{R}(R/\mathfrak{m})$}; \\
  \node(3) {$R/\mathfrak{m}$}; & \node (4) {$R/\mathfrak{m}$};\\};
  \draw[decoration={markings,mark=at position 1 with {\arrow[scale=1.5]{>}}},postaction={decorate},shorten >=0.5pt] (1) -- (2) node[anchor=south] {$a$};
  \draw[decoration={markings,mark=at position 1 with {\arrow[scale=1.5]{>}}},postaction={decorate},shorten >=0.5pt] (1) -- (3) node[anchor=west] {$\pi_{R/\mathfrak{m}}$};
  \draw[decoration={markings,mark=at position 1 with {\arrow[scale=1.5]{>}}},postaction={decorate},shorten >=0.5pt] (2) -- (4) node[anchor=west] {$\pi_{R/\mathfrak{m}}$};
  \draw[decoration={markings,mark=at position 1 with {\arrow[scale=1.5]{>}}},postaction={decorate},shorten >=0.5pt] (3) -- (4) node[anchor=south] {$a$};
\end{tikzpicture}\]
implies that the $R$-homomorphism $a1^{R/\mathfrak{m}}:R/\mathfrak{m} \rightarrow R/\mathfrak{m}$ is surjective. But this map is zero, so we arrive at a contradiction.
Hence, $\mathfrak{m}=\mathfrak{n}$.

The converse is immediate.
\end{prf}


\end{document}